\input amstex

%This will run under Ams-Tex (without style files

%The following (from Ams-Tex) need to be defined 
%or redefined for plain Tex

%\plainfootnote   \,   \frac   
%\text   \pmb (= Poor man's bold)   \Cal (calligraphic font)
% \comment  \endcomment   \cases \endcases
% \matrix \endmatrix

%%%%%%font definitions

\font\rm=cmr10 \rm

\font\bf=cmb10
\font\Rm=cmr9 at 11pt
\rm
\font\it=cmsl9 at 10pt
 %[track+70] 
 at 7pt

\font\Rrm=cmr17 at 16pt
   \font\Rm=cmr12 at 11.5pt
%%%%%%%%%%

  %%% (P,x^n)
\long\def\Pf{\par\noindent {\it Proof.} }
\def\({\left(}
\def\){\right)}
\def\st{such that }
\def\qed{\hfill$\bullet$\vskip 4pt}

\def\brcs#1{\left\{ #1\right\}}

\def\iso{\cong}
\def\wrt{with respect to }
\def\:{\,:}

\def\ker{\text{ker\,}}

\def\R{\text{\bf R}}
\def\N{\text{\bf N}}
\def\Z{\text{\bf Z}}
\def\Q{\text{\bf Q}}

\def\Arrow #1;#2.{#1\:#2 \to }
%#1: #2 -> 

\def\Set#1#2{\brcs{#1 \left|\vphantom{#1 #2} \right.#2}}
%set of #1 such that #2

\def\Oh#1{{\pmb O}\(#1\)}
%%%bigOh notation 

%%%little oh notation

%product left and right
\def\Rrr#1,#2{{\Cal J}_{#1,#2}}%difference
\def\slfrac#1#2{{\raise -.07 ex\hbox{$^{#1}$}}\!/\raise .35 ex \hbox{${}_{#2}$}}
\def\ssf #1/#2{\slfrac {#1}{#2}}

\def\pd #1,#2.{\frac {\partial #1}{\partial #2}}

%%%%%%%statement of lemmas, propositions,
%%theorems, examples, ...  \Lem 
   \long\def\Lem
#1.#2\par{\vskip4pt{\baselineskip=13pt\font\it=cmsl12 at
11.5pt\Rm
   \noindent {\rm \uppercase{#1}} #2\vskip3pt

   }}

\long\def\Title #1\par {\noindent{\Rrm #1}\vskip 9pt}
%main title

 \long\def\SubT #1.{\noindent {\it #1\/} } 
%subtitle
 
 \long\def\SecT
#1\par{\vskip 3pt \noindent {\bf #1}\vglue1pt
   \noindent}%section title

\long\def\subtitle #1.{\vskip 2pt \noindent {\it #1}}

\long\def\Rmk#1\par{\vskip 1pt \noindent {\it
Remark.} #1\vskip2pt}%remark

%%%%%%%%%macros for numbering
%%of propositions, lemmas,examples, theorems

\scrollmode\NoBlackBoxes
\magnification=1100
\long\def\Abstract #1\par%
{\vskip .2 true cm{\leftskip 1 true in \rightskip 1 true in \font\rm=cmr8 \rm
\baselineskip=1pt \font\it=cmsl8 \font\bf=cmb8
\parindent=0em {\bf Abstract} #1

}}
\comment
\font\rm=Times at 10pt

\font\bf=TimesB
\font\Rm=Times at 11pt
\rm
\font\it=TimesI at 10pt
 %!![track+70]
\endcomment

  %%% (P,x^n)
\long\def\Pf{\par\noindent {\it Proof.} }
\def\({\left(}
\def\){\right)}
\def\st{such that }
\def\qed{\hfill$\bullet$\vskip 4pt}

\def\brcs#1{\left\{ #1\right\}}
\def\Set#1#2{\brcs{#1 \left|\vphantom{#1 #2} \right.#2}}

\def\iso{\cong}
\def\wrt{with respect to }
\def\:{\,:}
%\long\def\Lem #1 #2\par{\noindent {\Sc lemma #1.} {\Rm
%#2}\vskip 2pt}
\def\Arrow #1;#2.{#1\:#2 \to }

\def\Oh#1{{\pmb O}\(#1\)}
%%%bigOh notation 

%%%little oh notation
\def\R{\text{\bf R}}
\def\N{\text{\bf N}}
\def\Z{\text{\bf Z}}
\def\Q{\text{\bf Q}}
 
%product left and right
\def\Rrr#1,#2{{\Cal J}_{#1,#2}}%difference

\def\slfrac#1#2{{\raise -.07 ex\hbox{$^{#1}$}}\!/\raise .35 ex \hbox{${}_{#2}$}}
\def\ssf #1/#2{\slfrac {#1}{#2}}

\def\pd #1,#2.{\frac {\partial #1}{\partial #2}}

    %notation: R_P, (P,x), pure traces etc
    %\C complexes \R reals \quotes#1
    %P will be of rad of cvg =1 with no negative coef 
%subordinate to: subeqv; f subeqv P means 
% all coefficients of P are nonnegative, and the 
%absolute values of coef of f are bdd above by etc.
    %equivalent to: eqv
    %explanation of 0 \leq f

   \long\def\Title #1\par {\noindent{\Rrm #1}\vskip 9pt}
 \long\def\SubT #1.{\noindent {\it #1\/} }   \long\def\SecT
#1\par{\vskip 3pt \noindent {\bf #1}\vglue1pt
   \noindent}%section title
\long\def\subtitle #1.{\vskip 2pt \noindent {\it #1}}

\long\def\Rmk#1\par{\vskip 1pt \noindent {\it
Remark.} #1\vskip2pt}

%%%%%%%%%macros for numbering

%5A
%7
%8
%9
%8A

\def\Z{\text{\bf Z}}

\def\ker{\text{ker\,}}

%this is the difference
%\def\rr#1{{\Cal V}\(#1\)}

\def\rank{\text{rank\,}}

\def\Z{\text{\bf Z}}

\def\ker{\text{ker\,}}

%this is the difference
%\def\rr#1{{\Cal V}\(#1\)}
\long\def\subsubtitle #1.{\vskip 3pt \noindent {\it #1 }}
\def\Aff{\text{Aff\,}}

\long\def\subtitle #1\par#2{\vskip 3pt \noindent {\it #1}\par \noindent
#2}

\def\ideal(#1){\(#1\)}
\def\oid(#1){< #1 >}

\Title Real dimension groups%
\plainfootnote{\rm${}^0$}{\rm Working document.}

\Abstract The main result is that dimension groups (not countable) that are also real ordered vector spaces can be obtained as  direct limits (over  directed sets) of simplicial real vector spaces (finite dimensional vector spaces with the coordinatewise ordering), but the directed set is not as interesting as one would like---e.g., it is not true that  a countable-dimensional real vector space which has interpolation can be represented as such a direct limit over the positive integers. It turns out this is the case when the group is additionally simple, and it is shown that the latter have an ordered tensor product decomposition. In the Appendix, we provide a huge class of polynomial rings that with a pointwise ordering are shown to satisfy interpolation, extending a result outlined by Fuchs. 

\vskip4pt 
\noindent {{\it David Handelman}}\plainfootnote{$^1$}{Supported in part
by a Discovery Grant from NSERC.}\vskip 4pt 

\noindent
Let $F$ be a subfield of the reals (although the real case of interest  occurs when $F = \R$), equipped with the relative ordering; $F^+$ will denotes the set of positive real numbers in $F$. A {\it partially ordered $F$-vector space\/} will be a vector space, $V$, over $F$, together with its positive cone, $V^+$, which satisfies the following properties:
$$V^+ + V^+ \subseteq V^+; V^+ - V^+ = V; V^+ \cap -V^+ = \brcs{0}; V^+ \cdot F^+ \subseteq V^+.
$$

We say a partially ordered $F$-vector space is {\it $F$-simplicial\/} (or simply {\it simplicial\/}) if there exists an integer $n$ \st it is isomorphic (as ordered $F$-vector spaces) to $F^n$ (the space of columns of size $n$ with the coordinatewise ordering acquired from $F$). Finally, a partially ordered abelian group (of which partially ordered $F$-vector spaces are special cases) $G$ satisfies {\it interpolation\/} if for all pairs of pairs $x_1, x_2$ and $y_1, y_2$ of elements of $G$ \st $x_i \leq y_j$ for all $i,j = 1,2$, there exists $z$ in $G$ \st $x_i \leq z \leq y_j$ for all $i,j$.

We wish to emulate  [EHS; Theorem 2.2]---that an unperforated partially ordered abelian group with interpolation is a direct limit of simplicial ordered groups (i.e., $\Z^n$ with the usual ordering, and positive group homomorphisms between them). In other words, we wish to characterize the direct limits of simplicial $F$-vector spaces (all direct limits are over directed sets, and with positive $F$-linear maps). Note that such objects are already dimension groups (direct limits of simplicial ordered groups), since unperforation is automatic for vector spaces over an ordered field (this requires inverses of nonzero positive elements be positive, which is trivially the case here). It turns out that there {\it is\/} a theorem, but it is not quite what it should be.

The following should come as no surprise.

\Lem Theorem 1. Let $V$ be a partially ordered $F$-vector space. Then $V$
 can be written as  a direct limit of $F$-simplicial partially ordered vector spaces if and only if $V$ satisfies interpolation.

Suppose however, that $V$ satisfies the hypotheses, and in addition, has countable dimension as an $F$-vector space. We expect a direct limit representation over a countable directed set, which is equivalent to a direct limit of the form
$$ F^{n(1) } @>M(1)>> F^{n(2)} @>M(2)>> F^{n(3)} @>M(3)>>\dots \tag 1
$$
where the maps ($M(i)$, $i=1,2,3, \dots$) between the simplicial vector spaces  are implemented by matrices all of whose entries are in $F^+$ (that is, the directed set can be taken to be the positive integers). In fact, typical representation theorems of   EHS type first deal with this countable version, and later extend to arbitrary direct limits when cardinality conditions are lifted. To obtain such a direct limit in the context of partially ordered vector spaces  requires an additional hypothesis.

We say a partially ordered $F$-vector space $V$ is {\it countably $F^+$-generated\/} if there exists a countable subset, $S$, of $V^+$ \st every element of $V^+$ is in the $F^+$-span of $S$---in other words, as an $F^+$-semigroup, $V^+$ is countably generated. It is not true that a countable-dimensional partially ordered $\R$-vector space is countably $\R^+$-generated, even if interpolation is thrown in (Example 7)! There may be similar cardinality problems for bigger dimensions, especially if  the continuum hypothesis is negated.

If $V$ is a direct limit as in (1), then it obviously satisfies countable $F^+$-generation.

\Lem Theorem 2. Let $V$ be a countable dimensional partially ordered $F$-vector space. Then $V$ is a direct limit over a countable index set of $F$-simplicial vector spaces if and only if $V$ is countably $F^+$-generated and satisfies interpolation.

Obvious examples of partially ordered $F$-vector spaces that are direct limits of simplicial ones include those of the form $G \otimes_{\Z} F$ where $G$ is a dimension group, and we take the ordered tensor product, as in [GH]. However, in these cases, the maps between the simplicial vector spaces are matrices all of whose entries are nonnegative integers (and this is a useless characterization of such limits). We later show that simple partially ordered $F$-vector spaces admit such a decomposition, giving an alternative proof of Theorem 1 for this class. 

The method of proof of Theorem 1 of course goes back to Shen's idea, which is rather simple to adapt here. However, because of the countability problem, we have a difficulty, which fortunately can be circumvented by the method in [G, Theorem 3.17].

\Lem Lemma 3. Let $V$ be a partially ordered $F$-vector space, let $\Arrow g;V_0. V$ be an $F$-linear positive map from a simplicial $F$-vector space $V_0$, and let $a$ be an element of the kernel of $g$. Then there exist a simplicial $F$-vector space $V_1$ together with  $F$-linear positive maps $\Arrow g_{01}; V_0.V_1  $ and $\Arrow h; V_1.V$ \st $g = hg_{01}$ and $g_{01}(a) = 0$.

\Pf Let $\brcs{e_i}_{i=1}^n$ be the standard basis for $V_0 \iso F^n$ with the usual ordering, write $a = \sum p_i e_i$ with $p_i $ in $F$. Chop the index set $\brcs{1,2,3, \dots, n}$ into three pieces, $S = \Set{i}{p_i > 0}$, $T = \Set{i}{p_i < 0}$, and $U = \Set{i}{p_i = 0}$. If either $S$ or $T$ is empty, it is clear what to do, and so we assume both are nonempty.

Let $e_i \mapsto x_i$ in $V^+$, so that $\sum_{i \in S} p_i x_i = \sum_{i \in T} |p_i|x_i$. Define $E_i = |p_i| e_i$ if $i \in S \cup T$ and $e_i$ otherwise, so the map $V_0$ is given by $E_i \mapsto |p_i| x_i$ if $p_i \neq 0$ and $E_i \mapsto x_i$ otherwise. Define the elements $y_i = |p_i| x_i$ in $V^+$ for all $i$. We now have the equation
$$
\sum_{S} y_s = \sum_{T} y_t.
$$
Since Riesz interpolation is equivalent to Riesz decomposition, there exist $y_{st}$ (with $(s,t)$ running over $S \times T$) in $V^+$ \st
$$\eqalign{
\text{for all $s \in S$,}\quad & y_s= \sum_{t \in T} y_{st} \qquad \text{and}\cr
\text{for all $t \in T$,}\quad & y_t = \sum_{s \in S} y_{st}.
}$$

Let $V_1$ be the simplicial $F$-vector space with standard basis given by $\brcs{f_{st}}_{S\times T} \cup \brcs{f_u}_{u \in U}$, and consider the assignment
$$\eqalign{
E_s &\mapsto \sum_T f_{st}\cr
E_t & \mapsto \sum_S f_{st}\cr
E_u& \mapsto f_u \cr
}$$

This extends to an $F$-linear positive homomorphism $\Arrow g_{01};V_0 . V_1$, since it is extendible from $e_i \mapsto \sum_j f_{ij}/|p_i|$ when $p_i \neq 0$.  Define $\Arrow h;V_1.V$ via $f_{st} \mapsto  y_{st}$ and $f_u \mapsto y_u$.

For $s \in S$,  $hg_{01} (e_s) = \sum_T y_{st}/p_s $, and this is $y_s/p_s = x_s = g(e_s)$; similarly $hg_{01}(e_t) = g(e_t)$ for $t \in T$; finally, for $u $ in $U$, $hg_{01} (e_u) = h(f_u) = y_u = g(e_u)$. Hence $hg_{01}  = g$. Next,
$$\eqalign{
g_{01}(a) &= g_{01}\(\sum_S p_s e_s \) - g_{01}\(\sum_T |p_t| e_t\) \cr
& = \sum_S p_s \sum_T f_{st}/p_s - \sum_T |p_t| \sum_S f_{st}/|p_t| \cr
& = \sum_S \sum_T f_{st} - \sum_T \sum_S f_{st} \cr
& = 0.
}$$
\qed

\Lem Lemma 4.  Let $V$ be a partially ordered $F$-vector space, let $\Arrow g;V_0. V$ be an $F$-linear positive map from a simplicial $F$-vector space $V_0$. Then there exist a simplicial $F$-vector space $V^1$ together with  $F$-linear positive maps $\Arrow g^{01}; V_0.V_1  $ and $\Arrow h; V^1.V$ \st $g = hg^{01}$ and $\text{ker} g \subseteq \text{ker} g^{01}  $.

\Pf Since $\text{ker} g$ is finite dimensional, we may iterate the previous construction.
\qed

\noindent {\it Proof of Theorems 1 and 2.} Now Theorem 2 follows in the usual way (index the generating set by the positive integers, etc). However, to prove Theorem 1 (results of this type are easiest to prove from the countable case, but as discussed earlier, this cannot be done here), we rely on the method of proof of [G;  Theorem 3.17]. Let $Y$ be the set of all finite subsets of $V^+$; this is a directed set, and the method of op cit (suitably modified by replacing $\Z$ by $F$, combined with Lemma 4 here shows that for we may construct simplicial $F$-vector spaces, $G_A$, together with positive $F$-linear maps $G_A \to G_B$ whenever $A \subseteq B$ that are compatible with the directed structure of $Y$, and positive $F$-linear maps $G_A \to V$ \st the image of  $G_A^+ $ in $V^+$ contains $A$, and such that any element of the kernel of $G_A \to V$ gets sent to zero in some $G_B$ under the map $G_A \to G_B$. As usual this means the limit over $Y$ is isomorphic to $V$, as partially ordered $F$-vector spaces. (Theorem 2 can also be proved from Theorem 1 by limiting $Y$ to finite subsets of a countable generating set, and observing that this is itself countable and directed, hence order final to the positive integers.)\qed

A stronger property than merely being a direct limit in the category of simplicial $F$-vector spaces, is that $V$ factorize as $W \otimes_{\Z} F$, where $W$ is a dimension group and  we take the tensor product ordering. Since direct limits factor through tensor products, such a decomposition automatically implies it is a direct limit of simplicial $F$-vector spaces. 
Moreover, it is routine that $\dim_F (W \otimes_{\Z} F) = \rank_{\Z} W$ (this has nothing to do with the ordering, and  reduces to the case that $W$ be a finite dimensional $\Q$-vector space, for which it is trivial).
In particular, if $V$ admits this factorization and is countable dimensional, then the corresponding $W$ must be countable, so that $V$ can be expressed as a direct limit with index set $\N$.

\Lem Proposition 5. If $V$ is a 
simple partially ordered $F$-vector space with interpolation, then there exists a dimension group $W$ \st $V$ is order-isomorphic (as partially ordered $F$-vector spaces) to $W \otimes_{\Z} F$.

\Pf  Pick a basis over $F$ for $V$, say $B:= \brcs{x_{\alpha}}$ \st at least one of the basis elements  is an order unit for $V$, and set $W $ to be the rational vector space spanned by $B$, with the relative ordering inherited from $V$. We note that even though traces are defined only as positive group homomorphisms to the reals, they are automatically $F$-linear (on an ordered $F$-vector space); this merely uses the fact that the rationals are order dense in $F$. 

The trace spaces of $W$ and $V$ are identical, that is, the inclusion induces an affine homeomorphism on the trace spaces---restriction is one to one since a trace is determined by its effect on $W$,  and conversely, every trace on $W$ extends to a trace on $V$, since $W$ contains an order unit of $V$.  

Since $W$ is a rational vector space, it is unperforated.  Suppose that $w$ is a nonzero element of $W^+ = W \cap V^+$; then $t(w) > 0$ for all traces of $V$, hence of $W$, and thus $w$ is an order unit. This yields that $W$ is simple, and moreover, that $W^+ \setminus \brcs{0}$ consists of the elements that are strictly positive at all traces of $V$. The trace space of $V$ is a Choquet simplex $K = S(V,u)$ where $u$ is a fixed order unit in $W$, and we have the representation  $W \to \Aff K$ obtained by restricting that of $V$.

The range of $V$ is dense in $\Aff K$, since $V$ is a simple dimension group unequal to $\Z$. Since $W$ is a rational vector space, so is its image in $ \Aff K$, and obviously the latter's closure contains the $F$-span of the image of $W$, in particular, the image of $V$. Since the latter is dense, so is the image of  $W$ in $\Aff K$. It follows that  $W$ is a dimension group. Next, consider the map $W \otimes F \to V$ (given by $\sum w_r \otimes_{\Z} \lambda_r \mapsto \sum w_r \lambda_r$). This is obviously an isomorphism of $F$-vector spaces (because of the construction), but we have to show it is an order isomorphism.

The positive cone of the tensor product consists of sums of terms of the form $w_r \otimes \lambda_r$ where $w_r \in W^+$ and $\lambda_r \in F^+$. These clearly go to  positive elements in $V$. Since $F$ has unique trace, it is easy to see that $H = W \otimes F$ has the same traces as $W$ and thus of $V$, that is, determined by $w \otimes \lambda \mapsto t(w)\lambda$. So we have a vector space isomorphism $H \to V$ that is positive, and induces a homeomophism on the trace spaces; moreover, $H$ and $V$ are simple (and unperforated), and it is immediate that this must be an order isomorphism (if $h \in H$ goes to something positive and not zero, its image must be strictly positive under all traces; from the homeomorphism, this pulls back to $H$, hence $h$ is strictly positive at all traces, and thus is positive).
\qed

\Lem Corollary 6. Suppose $V$ is a countable dimenensional 
simple partially ordered $F$-vector space with interpolation.  Then $V$ is order isomorphic to a direct limit of the form (1).

This suggests the question as to whether every dimension group that is also an ordered $\R$-vector space can be so decomposed (as $W \otimes_{\Z} \R$). We will see below that a plausible countable dimensional real algebra that is a dimension group (and an ordered ring) cannot be, because its positive cone is not countably $\R^+$-generated

If we iterate tensoring with $F$ (say  $F = \R$), we create monstrously large vector spaces; for example $\R \otimes_{\Z} \R$ is a real vector space $V $ of dimension $2^{\aleph_0}$, with a nonzero linear map $\Arrow \tau; V.\R$ \st $V^+ \setminus \brcs{0} = \tau^{-1} (\R^{++})$; this comes from the fact that $\R$ is merely a vector space of dimension $2^{\aleph_0}$ over the rationals, that the ordered tensor product of two dimension groups is a dimension group, that the ordered tensor product of two simple dimension groups is a simple dimension group, and that the trace space of an ordered tensor product  of dimension group is the corresponding tensor product in the category of Choquet simplices, in this case both consist of a singleton, hence the tensor product has unique trace. In this case, there is a simple characterization of $\R \otimes \R$ an {\it ordered  group\/}:  a rational vector space of dimension $2^{\aleph_0}$ having unique trace \st the image of whose values is all of $\R$.

Now we investigate necessary conditions for countable (and other) $F^+$-generation. Suppose that $(V,u)$ is a partially ordered $F^+$-vector space with order unit, and let $\,\hat{} \,$ denote the natural homomorphism $\Arrow {\hat{}} ; V . \Aff S(V,u)$. For every element $g$ in $V^+$, define the {\it zero set\/} of $g$, $Z(g) = \Set{\tau \in \partial_e S(V,u)}{\tau(g) = 0}$. Each $Z(g)$ is of the form $\partial_e L = L \cap \partial_e S(V,u)$ for $L$ a closed face of $S(V,u)$. Note that if $V$ is simple, then the only zero sets are the trivial ones, $\emptyset$ and $\partial_e S(V,u)$. 

\Lem Lemma 7. Suppose the partially ordered $F$-vector space with order unit $(V,u)$ has the property that $V^+$ is $F^+$-generated by a set of infinite cardinality $\aleph$.{\par}
\item{(a)} The number of subsets  of $\partial_e S(V,u)$ of the form $Z(g)$ (with $g$ in $V^+$) is at most $\aleph$;
\item{(b)} the number of order ideals of $V$ having their own (relative) order unit is at most  $\aleph$;
\item{(c)} the number of order ideals of $V$ is at most $2^{\aleph}$.

\Pf (a) Obviously, if $g = \sum r_i g_i$ where $\brcs{r_i} $ is a finite set of positive real numbers and $g_i$ are in $V^+$, then $Z(g) = \cap_i Z(g_i)$. Hence the zero set of an arbitrary element of $V^+$ is a finite intersection of zero sets of the generating set, hence the number of zero sets (of elements of $V^+$) is bounded above by the number of finite subsets of a set with $\aleph$ elements. Since $\aleph$ is infinite, the number of its finite subsets is again $\aleph$. 

(b) If $H$ is an order ideal with order unit $h$, we can write $h = \sum g_i r_i$ where $\brcs{r_i}$ is a finite set of positive real numbers and $g_i$ are in $V^+$. Since $g_i \leq h/r_i$, it follows that $g_i$ all belong to $H$, and obviously determine $H$ (the order ideal generated by the finite set $\brcs{g_i}$ is $H$). This yields an onto map from the finite subsets of the generating set to the set of order ideals with  order unit, hence the result.

(c) Trivial, and useless anyway.
\qed

\Lem Example 8. A countable-dimensional real ordered vector space which is a dimension group that cannot be represented as a sequential direct limit of simplicial $\R$-vector spaces. 

\Pf  Let $R = \R[x]$, the polynomial ring with the coordinatewise ordering on the unit interval. This is countable dimensional, but cannot have a countable $\R^+$-generating set---if it did, then with $\aleph = \aleph_0$ in Lemma 7, there would only be countably many zero sets of positive elements. However, the pure trace space consists of the point evaluations from points in the unit interval, and each element of $\brcs{(x-t)^2}_{t \in [0,1]}$ yields a zero set consisting of the singleton, the point evalution at $r$; so uncountably many zero sets exist. 

Now  $R$ has interpolation. This has an interesting history. It is stated without proof or references in [F1] and [F2], but in [F3, pp\,19--20] (I am indebted to George Elliott for finding this reference), a proof is sketched. Because the last is not readily available (even on-line), we include  Fuchs' result, and extend it to cover subfields of the reals in the Appendix (for the reals themselves, the proof simplifies to what Fuchs' had outlined).  \qed

\vskip 2pt \noindent {\bf Appendix} Interpolation for polynomial rings with the pointwise ordering

\vskip 4pt \noindent Let $I =[\gamma,\delta]$ be a closed bounded interval with interior, in $\R$. Let $L$ be a subfield of the reals, and form $R = L[x]$, the polynomial ring. Define a positive cone on the latter to be the set of elements of $R$ that are nonnegative on $I$. This of course depends on $L$, $\gamma,$ and $\delta$, so if there is any possible ambiguity, we write $R_{L, \alpha, \beta}$. If $L = \R$, then the rings are order isomorphic (via $x \mapsto cx + d$, $c \neq 0$) regardless of the choice of $\gamma < \delta$. At the other extreme, if $L= \Q$, then there are uncountably many order-isomorphism classes, as the only ring automorphisms that implement an order isomorphism are given by  via $x \mapsto cx + d$, $c \neq 0$, with $c,d$ rational (we can recover some of the properties by ordered ring invariants, e.g., the number of rationals in the set $\brcs{\gamma,\delta}$---that is, $0$, $1$, or $2$---is equal to the number of codimension 1 order ideals in $R_{\Q,\gamma,\delta}$).

We will show that for all choices of $L$, $\gamma < \delta$, the resulting ordered ring satisfies Riesz interpolation. If $L = \R$, Fuchs has given a proof, outlined in [F3, pp19--20], which uses Hermite interpolation, but really only requires the Chinese remainder theorem. On the other hand, if $L = \Q$, then there are many more technicalities (concerning behaviour of algebraic conjugates, etc). We use two key steps from Fuchs' argument.  If one wants a proof just for the reals, the following simplifies considerably to Fuchs' outlined proof.

Recall the Chinese remainder theorem: if $\brcs{I_s}_S$ is a finite collection of ideals in a ring $R$ \st  for all $s\neq t \in S$, $I_s + I_t  = R$, then the natural map $R \to \prod_S R/I_s$ is onto (usually stated in the equivalent form, $R/\cap_S I \to \prod_S R/I_s $ is an isomorphism).

Suppose $f_i \leq g_j$ (as functions on $I$), $i,j = 1,2$ with all four polynomials in $R$; to find $w$ in $\R$ \st $f_i \leq w \leq g_j$, we may assume $f_2 $ is identically zero (by subtracting $f_2$ from the other terms). We may also assume that the four polynomials are distinct. In this case, the set of zeros of the equation $f_1 \vee 0 = g_1 \wedge g_2$ is finite; let $S$ denote the set of these zeros that lie in $I$. We note that $S$ consists of real numbers that are algebraic over $L$.

If $\alpha$ and $\beta$ are roots of the same irreducible polynomial over $L$, we say $\alpha$ and $\beta$ are (algebraic) {\it conjugates\/} (over $L$); when both are sitting inside the reals (as is the case with the elements of $S$), then there is a field isomorphism $\Arrow \sigma; K:= L[\alpha]. L[\beta] \subseteq \R$ fixing $L$ pointwise, induced by $\alpha \mapsto \beta$. In the case that $L = \R$, this can be ignored.

Let $T$ be a cross-section of the conjugacy classes of elements of $S$---in other words, $T$ is a subset of $S$ \st each element of $S$ is conjugate to exactly one element of $T$; for technical reasons, we also choose  elements of $T$ to be in the interior of $I$ if the conjugacy class contains such an element. Obviously, if $L = \R$, or more generally, if all the zeros lie in $L$, then $S = T$.
For each $t \in T$, define the class $S_t = \Set{s \in S}{s \text{ is conjugate to } t}$, so that $\brcs{S_t} $ is a partition of $S$.  

The proof is divided in three.

\noindent Step I. {\it For each $t$ in $T$, there exists $h_t$ in $L[x]$ together with a relatively open neighbourhood in $I$, $U_t$, of $S_t$, \st $0,f_1 \leq h_t \leq g_1, g_2$ as functions on $I$ restricted to $U_t$.} This is a somewhat technical (but not difficult) result (routine over  $\R$), which will take the most time.

\noindent Step II. {\it There exists $h$ in $L[x]$ together with a relatively open neighbourhood in $I$, $U$, of $S$, \st $0,f_1 \leq h \leq g_1, g_2$ as functions on $U$.} This is the first step of Fuchs' argument over the reals (obtained by appealing to Hermite interpolation), and here it follows from Step I by the  Chinese remainder theorem.

\noindent Step III. {\it There exists a  polynomial, nonnegative on $I$, $f$ in $L[x]$, with zeros on a finite set including  $S$ \st $-h/f,(f_1-h/f) \leq (g_i - h)/f$ on all of $I$, except  the zeros of $f$; additionaly, these inequalities   can be interpolated by an element of $L[x]$.} This is virtually identical to Fuchs' second step, and easily concludes the proof.

We now proceed to the proof of Step I. For each $t \in T$, let $p_t$ denote the unique monic polynomial in $L[x]$ irreducible over $L$ satisfied by $t$; automatically, it is also the irreducible polynomial of all  $s$ in $S_t$. Now we make an observation about local nonnegativity of polynomials at algebraic numbers.

Suppose $\alpha $ is algebraic over $L$ with monic irreducible polynomial $p$ in $L[x]$, and $\alpha \in I$. Let $c$ in $L[x]$ be nonconstant and satisfy $c (\alpha) = 0$. Then there is a relative neighbourhood, $U$, of $\alpha$, in $I$ (that is, we have to include relatively open sets of the form $[\gamma, \gamma + \epsilon)$ and $(\delta-\epsilon, \delta]$) \st $c|U \geq 0$ if and only if there exists an integer $k$ and an element $N$ of $L[x]$ with $N(\alpha) \neq 0$ \st $c = p^k N$ \st
$$
\cases
k \text{ is even and $N(\alpha) > 0$}& \text{if $\alpha \in (\gamma,\delta)$} \\
N(\alpha) > 0 & \text{if $k $ is even and $\alpha \in \brcs{\gamma,\delta}$}\\
(p'(\alpha))^k N(\alpha) > 0& \text{if $k $ is odd and $\alpha= \gamma$}\\
(p'(\alpha))^k N(\alpha)< 0& \text{if $k $ is odd and $\alpha= \delta$}.\\
\endcases
$$
This follows easily from  $p(\alpha) = 0$ implies $D^k (p^k N)(\alpha) = k! (p'(\alpha))^k N(\alpha)$. In the latter two cases ($k$ is odd, and $\alpha$ is an endpoint), we may obviously replace $(p'(\alpha))^k N(\alpha)$ by $p'(\alpha) N(\alpha)$. Of course, we really don't need $k \geq 1$, it obviously also applies with $c(\alpha) \neq 0$.

\Lem (Technical) Lemma. For $t$ in $T$, there exist an integer $m \equiv m(t)$ to together with a family of nonempty  open intervals in the reals, $\brcs{(a(s), b(s))}_{s\in S_t}$ and  $\brcs{r_j}_{j=1}^{m-1}$ in the extension field $L[t]$ \st if $h$ is any  polynomial with coefficients from $L$ satisfying
  {\par}
\item{(i)} $h(t) =\inf \brcs{g_1, g_2} (t)$,
\item{(ii)} $h^{(j)}(t) = r_j$ for $1 \leq j \leq m-1$,
\item{(iii)} for all $s$ in $S_t$, $h^{(m)}(s)  \in (a(s), b(s))$,
{\par\par}
\noindent then there exists a relative neighbourhood $U_t \equiv U(t,h)$ of $S_t$ \st as functions on $U_t$, $f_1,0 \leq h \leq g_j$.

\Pf Case 1. $f_1 (t) > 0$. Then $\inf \brcs{g_1 (t),g_2(t)} = f_1(t) > 0$; in particular, $f_1 (s) \neq 0$ for any algebraic conjugate $s$ (i.e., for $s \in S_t$); from $g_i(t) = f_1(t)$ (which holds for at least one of the $g_1, g_2$),  we have $f_1(s) = g_i(s) \neq 0 $; as $g_i \geq 0$, $g_i (s) >0 $, so that $f_1 (s) > 0$, whence $g_{2-i} (s) > 0$. In particular, for all $s$ in $S_t$, we have $\inf \brcs{g_1 (s),g_2(s)} = f_1(s) > 0$, and obviously $g_i (t) = f_1 (t)$ if and only if for all $s$ in $S_t$, if and only for some $s \in S_t$, we have $g_i (s) = f_1 (s)$. In particular, by (i), $h(s) = f_1(s) > 0$ for all $s$ in $S_t$, so there exists a neighbourhood of $S_t$ on which $h$ is nonnegative. 

\vskip2pt \noindent Case 1a. {\it Assume $t \in (\gamma,\delta)$.} Let $p$ be the monic irreducible polynomial with coefficients from $L$ satisfied by $t$ (hence by all other elements of $S_t$).  From $g_i - f_1 \geq 0$, at least one vanishing at $t$, we have a factorization $g_i - f_1=p^{2m_i} M_i$ where $m_i$ are nonnegative integers with at least one being nonzero, and $M_i$ are in $L[x]$ with $M_i(t) > 0$. Since $g_i - f_1 \geq 0$ on all of $I$, this forces $M_i (s) \geq 0$ for all $s$ in $S_t$; since $M_i (t) \neq 0$, we have $M_i (s) \neq 0$ for every such $s$, and therefore $M_i (s) > 0$ for all $s \in S_t$.

Set $m = \max\brcs{2m_1,2m_2} $ and $r_j = D^j(f_1) (t)$. Without specifying the values of $h^{(m)(s)}$ yet, we see that $h-f_1$ has a zero of order (at least) $m$ at $t$, and thus we can write it as  $p^m P$ where $P$ is a polynomial (obviously depending on our choice of $h$), and $h^{(m)} (s) =m! (p'(s))^m P(s) + f_1^{(m)}(s)$ for all $s \in S_t$) (since all fields here are separable---!---it follows that $p'(s) \neq 0$); since $m$ is even, we see that  the condition that $D^m (h-f_1) (s) > 0$ is sufficient for $h \geq f_1$ on a neighbourhood of $s$; this translates to $h^{(m)}(s) > f_1^{(m)} (s)$. 
Thus  we let $a(s) = f_1^{(m)}(s)$. The set of  conditions, $h^{(m)}(s) > f_1^{(m)} (s)$ (one for each $s$ in $S_t$) is sufficient to guarantee $0, f_1 \leq h$ on a neighbourhood of each $s$ in $S_t$. (Notice that there is no problem if one or both of the endpoints are included in $S_t$.)

Now consider $g_i - h = (g_i - f_1) +( f_1 - h) = p^{2m_i} M_i - p^m P$; if for one of the $i $, $m_i < m/2$, then dividing by $p^{2m_i}$ we see that $M_i (s) > 0$ is sufficient for $g_i - h \geq 0$ on a neighbourhood of $s$, and this automatically follows from the first paragraph for every $s$ in $S_t$.
In addition, in this case,  we must have $m_{2-i} = m/2$, so that $g_{2-i} - h = p^m(M_i - P)$; thus  sufficient for $g_{2-i} - h$ to be nonnegative on a neighbourhood of $s$ (again, since $m$ is even) is that $M_{2-i} (s) - P(s) > 0$. Since $P(s) = (h^{m}(s)- f_1^{(m)}(s))/p'(s)^m m!$, so we just need to ensure that $h^{(m)} (s) < m! p'(s)^m M_{2-i}(s) + f_1^{(m)}(s)$ (since $m$ is even, the sign of $p'(s)$ is irrelevant). In this case, we set $b(s) = m! p'(s)^m M_{2-i}(s) + f_1^{(m)}(s)$ (where $i$ is defined by $m_i < m/2$). Now we must check that $a(s) < b(s)$, that is, $0 < m! p'(s)^m M_{2-i} (s)$, which is obvious since (again) $m$ is even.

There remains the possibility that  $m_1 = m_2 = m/2$; and the same analysis yields a choice for $b(s)$, namely $m! p'(s)^m \min\brcs{M_1(s), M_2(s)} + f_1^{m} (s)$.

\vskip2pt\noindent Case 1b. {\it $t = \gamma$ or $\delta$ but $|S_t| = 1$.} The process for $t = \delta$ is obtained from the process for $t = \gamma$ by applying an automorphism of $L[x]$ which reverses the orientation (possibly shifting the interval at the same time, and we find that the definitions of $a(s)$ and $b(s)$ are obtained in reverse to the way they were obtained  in the other subcases), so we can just assume that $t = \gamma$ and $t$ has no conjugates in the interval $(\gamma,\delta]$. We write $g_i - f_1 = p^{k_i} M_i$ with $M_i (t) > 0$ if either $k_i$ is even or if both $k_i $ is odd and $p'(t) >0$, and $M_i(t) < 0$ if both $k_i$ is odd and $p'(t)< 0$. To simplify matters, we replace $p$ by $-p$ if $p'(t)< 0$, that is, we can assume $p'(t) > 0$ if we don't mind losing monicity, and eliminating the last possibility---so that $M_i (t) > 0$ in any case. Set $m = \max\brcs{k_1, k_2}$, and $r_j = f_1^{(j)}(t)$. Then $h$ satisfying (i) and (ii) entails $h - f_1$ has order at least $m$ at $t$, and thus factors as $p^m  P$; as in Case 1a, we see quickly that $h^{(m)}(t) > f_1^{(m)}(t)$ is sufficient for $ h \geq f_1, 0$ on a relative neighbourhood of $\gamma$ (that is, of the form $[\gamma,\gamma + \epsilon)$). 

As in Case 1a, we can write $g_i - h = p^{k_i} M_i - p^m P$; now the fact that $M_i (t) > 0$ and $p'(t) > 0$ allows the same analysis as in Case 1a, to show that we may choose $b(s) = m! p'(s)^m M_{2-i}(s) + f_1^{(m)}(s)$ if $k_i < k_{2-i} = m$, and $b(s) = m! p'(s)^m \min\brcs{M_1(s), M_2(s)} + f_1^{m} (s)$ if $k_1 = k_2  = m$.

\vskip2pt\noindent Case 1c. {\it $S_t = \brcs{\gamma,\delta}.$} In this weird case, we have (labelling $t = \gamma $
and $s = \delta$), $g_i - f_1 \geq 0$ entails $g_i - f_1 = p^{k_i} M_i$ where if $k_i$ is even, then $M_i (s), M_i(t) > 0$, while if  $k_i$ is odd,  then $M_i (t) p'(t) > 0$ and  $M_i(s) p'(s) < 0$. Let $m = \max\brcs{k_1, k_2}$ and again set $r_j = f_1^{(m)}(t)$, so that $h - f_1 = p^m P$. 

Suppose $m$ is even; then sufficient for nonnegativity of $h-f_1$ on a neighbourhood of $S_t$ is that $P(s), P(t) > 0$. Since $h_1^{(m)}(\alpha) - f_1^{(m)}(\alpha) = m! p'(\alpha))^m P(\alpha)$ for $\alpha \in S_t$, as in all the previous subcases, we can set $a(t) = f_1^{(m)}(t) $ and $a(s) = f_1^{(m)}(s)$. Continuing with even $m$, the same arguments as in the previous subcases give the same choice for $b(t)$ and  $b(s) $.

Finally (for Case 1), suppose $m$ is odd. We want to ensure $p'(t)P(t) > 0$ and $p'(s) P(s) < 0$. Since $S_t$ consists only of $t$ and $s$ and nothing in between, we see that they are consecutive real roots (each of multiplicity one) of the real polynomial $p$; hence $p'(t)p'(s) < 0$ (signs of the derivatives are opposite); this is convenient. By replacing $p$ by $-p$ if necessary, we may assume $p'(t) > 0$, and thus $p'(s) < 0$. 
%Hence we have to impose conditions on $h^{(m)}(s)$ and $h^{(m)}(s)$ to guarantee $P(s), P(t) > 0$.
 We want to ensure that $D^m (h-f_1)(t) > 0$ and $D^m (h-f_1)(s) < 0$. Set $a(t) = f_1^{(m)}(t)$ and $b(s) = f_1^{(m)}(s)$; note the appearance of $b(s)$ {\it not\/} $a(s)$. Now similar analysis with the $g_i$ as in all the previous subcases (I'm getting pretty tired) realizes the complementary $b(t) $ and $a(s)$.

There are no other subcases to consider, since we picked the cross-section $T$ so that if a conjugacy class contains an interior point of the interval, then we chose the corresponding representative in $T$ to be an interior point.

\vskip4pt\noindent Case 2. $f_1 (t) = 0$. There exists $i$ 
\st $g_i (t) = 0$, hence $g_i(s) = 0$ for all $s$ in $S_t$. 

\vskip2pt\noindent Case 2a. {$t \in (\gamma,\delta)$.} If $g_i (t) > 0$ (necessarily $g_{2-i} (t) = 0$), then $g_i(s) \neq 0$ for all $s$ in $S_t$, so from $g_i \geq 0$, we have $g_i (s) > g_{2-i} (s) = 0$, hence $g_i > g_{2-i}$ on a neighbourhood of $S_t$. Hence we can disregard $g_i$---we only have to guarantee that $g_{2-i} \geq h\geq f_1,0$ on a neighbourhood of $S_t$. If the order of $t$ as a zero of $f_1$ is odd, it must be at least as large as the order of $t$ for $g_{2-i}$, since $g_{2-i} - f_1 \geq 0$. Dividing by a sufficiently high (but nonzero) even power of $p$, $p^{2l}$ that divides both $g_{2-i}$ and $f_1$, we reduce to the situation that either $f_1/p^{2l} (s) <  0$ or $g_{2-i}/p^{2l} (s) > 0$; in the former situation, set $m = 2l$ and $r_j = 0$, and verification is trivial (we can take $b(s) = \infty$
for every $s$ for which $f_1/p^{2l} (s) <  0$), and in the latter case, if it occurs for one value of $s$ in $S_t$, it occurs for all, and then we either have $g_{2-i}/p^{2l} (t) = f_1/p^{2l} (t) > 0$, which is Case 1, or $g_{2-i}/p^{2l} (t) > f_1/p^{2l} (t) $, which isn't any case at all. In every single one of these possibilities, the choices for the intervals $(a(s), b(s)) $ are straightforward.

If  both $g_i (t) = 0$,  the order of $t$ as a zero of each $g_i$ is  even (as $t$ is an interior point), and of course $g_i - f_1 \geq 0$; if the order of $t$ as a zero of $f_1$  is odd, its order must be at least as large (and therefore more than  the order at $g_i$. We may thus divide everything  by $p^{2}$  (not affecting any of the inequalities, since $p^2$ is nonnegative), and continue this process as far as possible. At that point, either we reduce to Case 1, or to not both $g_i (t) $ being zero, the situation of the previous example, and again the choices for $r_j$ and $m$ are routine.

\vskip2pt\noindent Case 2b. {\it $t =\gamma$ or $\delta$ but $|S_t| = 1$.} Reduce to $t= \gamma$ as in Case 1b. Replace $p$ by $-p$ if necessary,  to ensure that $p'(\gamma) > 0$. Then $p^k \geq 0$ on $I$, since $p$ has no zeros on $(\gamma,\delta]$. Hence we may proceed as in Case 2a, not worrying about the parity of the power of $p$.

\vskip2pt\noindent Case 2c. {\it $S_t =\brcs{\gamma,\delta}$.} Again replace $p$ by $-p$ if necessary to ensure that $p'(\gamma) > 0$, so that $p$ is strictly positive on the interior of $I$, and thus $p|I \geq 0$, so we can again proceed as in Case 2a, dividing by a power (not worrying about the parity) to reduce to other cases.
\qed

\Lem Lemma. For $r$ a real number that is algebraic over  $L \subseteq \R$ of degree $n$, set $R = L[x]$ and let $K $ be the $n$-dimensional field extension $L[r]$. Let $m$ be a positive integer. The map
$\Arrow \phi;L[x] . K^{m+1}$ defined by $f \mapsto (f^{(j)}(r))_{j = 0}^m$, is onto.

\Pf Let $p$ be  a minimal polynomial of $r$ \wrt $L$, so that $\deg p = n$. Then the kernel of $\phi$ is exactly the ideal  $J = p^{m+1}R$ of $R$ and $\phi$ is $L$-linear; it thus induces the one to map  $L$-linear map $\Arrow \overline \phi; R/J . K^{m+1}$. The $L$-dimension of the left side is $n(m+1)$, so that the left side is equidimensional (as an $L$-vector space) with the right side, hence $\overline \phi$ is an isomorphism. \qed

The following observation, when applied to $L = \Q$, is well known. Let $K$ be a formally real finite dimensional extension field of $L$, and let $\brcs{\sigma}_{\Sigma}$ be a  family of homomorphisms $K \to \R$ (at least one exists by formal reality). Then the image of $K$ under the obvious map $K \to \R^{\Sigma}$ is dense. We note that $\Sigma$ is automatically linearly independent, so both $|\Sigma| <\infty$ and density follow immediately.

\Lem  Lemma. (Step I) For each $t$ in $T$, there exists $h_t$ in $L[x]$ together with a relatively open neighbourhood in $I$, $U_t$, of $S_t$, \st $0,f_1 \leq h_t \leq g_1, g_2$ as functions on $I$ restricted to $U_t$.

\Pf For each $s$ in $S_t$, choose field isomorphisms (over $L$) $\Arrow \sigma_s; L[t]. L[s] \subseteq \R$ \st $\sigma_s (t) = s$, and let $\Sigma$ be the set of such $\sigma_s$ (this will of course vary as $t$ varies). The map $L[t] \to \R^{\Sigma}$ has dense range, as we observed earlier. Hence we may find $q  $ in $L[t]$ \st for all $s$ in $S_t$, we have $\sigma_s j\in (a(s), b(s))$. With $K = L[t]$, apply the previous lemma to the sequence $(f_1(t), r_1, \dots, r_{m-1}, q$. Ontoness of the map entails that there exists $h_t$ in $L[x]$ satisfying the technical lemma, hence $0,f_1 \leq h_t \leq g_i$ on a relative neighbourhood of $S_t$ in $I$.
\qed

\Lem  Lemma. (Step II) There exists $h$ in $L[x]$ together with a relatively open neighbourhood in $I$, $U$, of $S$, \st $0,f_1 \leq h \leq g_1, g_2$ as functions on  $U$.

\Pf For each $t$ in $T$, let $p_t$ be an irreducible polynomial in $L[x]$ for $t$; the corresponding $m_t$ comes from the technical lemma above. Set $J_t = p_t^{m_t}$, so that the $\brcs{J_t}_{t\in T}$ are pairwise comaximal. By the Chinese remainder theorem, there exists $h$ in $L[x]$ \st for each $t$, $h - h_t \in J_t$. But this simply means that $D^{k} h(t)= D^k h_t (t)$ for $0 \leq k \leq m_t$, and since all the derivatives of $h$ and $h_t$ belong to $L[x]$, we also have $D^{m_t} h (s) = D^{m_t} h_t (s)$ for all $s$ in $S_t$. Hence $h$ satisfies the conditions of the technical lemma for each equivalence class, and thus $0 ,f_1\leq h \leq g_1,g_2$ on a neighbourhood of $S = \cup_{t\in T} S_t$. 
\qed 

\Lem Lemma. (Step III) There exists a nonzero polynomial, nonnegative on $I$, $f$ in $L[x]$, together with a finite subset $S_0$ of $I$ containing $S$ \st $-h/f,(f_1-h)/f \leq (g_i - h)/f$ on $I \setminus S_0$ and can be interpolated by an element of $L[x]$. 

\Pf Let $S_0$ be the union of the zero sets of  $h, f_1 -h, g_i-h$ intersected with $I$. Obviously $S \subseteq S_0$, and none of $S_0 \setminus S$ can be conjugate to an element of $T$. For each conjugacy class in $S_0$, with representative $u$, let $p_u$ be a minimal polynomial of $u$ over $L$; let $T_0$ be a cross-section of $S_0$ (enlarging $T$). Let $M$ be any even number exceeding the orders of all the zeros of the four polynomials, so that with $f =\prod_{T_0} p_u^M$, we have that $f|I \geq 0$ and moreover, the poles of {\it each\/} of $h/f, (f_1-h)/f , (g_1 - h)/f, (g_2 -h)/f$ in $I$ occur exactly at the points of $S_0$ and no others. 

In a (relative) neighbourhood of $s \in S_0$, we must have $\lim_{x\to s} (g_i - h)/f = +\infty$ (the limit is two-sided if $s$ is not a boundary point, but only one-sided if it is one of $\gamma$ or $\delta$) since $g_i - h $ is nonnegative on $I\setminus S$. Since $-h$ and $f_1 - h$ are negative on $I\setminus S_0$, it similarly follows that $\lim_{x\to s} - h/f = \lim_{x\to s} (f_1 - h)/f = -\infty$. On the other hand, we note that all of  $d_0 = \sup_{x \in I} -(h/f)(x)$, $d_1 = \sup_{x \in I} ((f_1-h)/f)(x)$, $e_i = \inf_{x \in I} ((g_i - h)/f)(x)$ are finite. There  exists a relative neighbourhood $V$ of $S_0$ \st on $V \setminus S_0$, all of $(g_i- h)/f >1 +  \max\brcs{d_0,d_1} $ and $-h/f, (f_1-h)/f < \inf\brcs{e_i} - 1 $ hold. On the remainder of $I$, $((g_i - h)/f) (x) > -(h/f)(x), (f_1-h)/f)(x)$; by compactness of $I\setminus V$, there exists $\eta > 0$ \st all four differences are at least $\eta$ on $I \setminus V$. Set $G_i = (g_i - h)/f \wedge (1 + \max \brcs{d_0,d_1})$, $E_1 = -h/f \vee (\min\brcs{e_i} - 1)$, and $E_2 = (f_1 -h)/f \vee (\min\brcs{e_i} - 1)$ so that $G_i, E_i $ are all continuous, and $G_i - E_j > \min\brcs{1, \eta}$ on all of $I$. Let $\kappa = \min\brcs{1, \eta}$, and set $G = G_1 \wedge G_2$ and $E = E_1 \vee E_2$. Then $G - E > \kappa$ on $I$

By the Weierstrass density theorem, there exists a {\it real\/} polynomial $w$ \st \wrt the sup norm, $ \| w - (G + E)/2 \| < \kappa/4$. Since $L$ is dense in $\R$, there exists a polynomial $W$ in $L[x]$ \st $\| W - w\| < \kappa /8$. As $\| G -(G + E)/2 \| \geq \kappa/2$, it easily follows that $G \geq W \geq E$ on $I$, and $W$ is thus the desired interpolant. \qed

\Lem Theorem. The ring $R_{L,\gamma,\delta} = L[x]$ equipped with the pointwise ordering from the interval $[\gamma,\delta]$ satisfies the Riesz interpolation property. 

\Pf From Steps I--III, we can find $W$ \st $-h/f , (f_1-h)/f \leq W \leq (g_i - h)/f$ on $I\setminus S$, with $f$ a square. Hence $0,f_1 \leq h +Wf \leq g_i$ on $I\setminus S$, and by continuity on all of $I$. \qed

The use of the Chinese remainder theorem together with the obvious result on values of derivatives in Step II yields a version of Hermite interpolation over the subfield $L$. We have to be careful here, because if for example, $\rho$ is transcendental over $L$ and $f$ is in $L[x]$, then prescribing a value for $f (\rho)$ either is impossible or  determines $f$ (and thus all of its values) completely! So we restrict the parameters to numbers algebraic over $L$ (and again noting that if $f^{(j)}(\alpha)$ is given, then $f^{(\beta)}$ is uniquely determined whenever $\beta$ is an algebraic conjugate of $\alpha$). The following unexciting result is proved as Step II was. 

\Lem  Proposition. (Hermite interpolation over real subfields) Let $L$ be a subfield of the reals. Suppose $\brcs{\alpha}_{\alpha \in Y}$ is a finite set of real numbers that are algebraic over $L$. Suppose that for each $\alpha$, there exists $m \equiv m(\alpha)$ together with $\brcs{r_{j,\alpha}}$, $0 \leq j \leq m(\alpha)$ with the following properties.
{\par}
\item{(i)} For each $\alpha$, all $r_{j,\alpha}$ belong to the extension field $K_{\alpha} := L[\alpha]$
\item{(ii)} If  $\alpha$ is conjugate to $\beta$, and $\Arrow \sigma;K_{\alpha}. K[\beta]$ is the field isomorphism induced by $\alpha \mapsto \beta$, then $m(\beta) = m(\alpha)$ and $r_{j,\beta} = \sigma(r_{j,\alpha})$ for all $0 \leq j \leq m(\alpha)$.
{\par}
\noindent Then there exists a polynomial $h$ in $L[x]$ \st for all $\alpha$ in $Y$, for all $0 \leq j \leq m(\alpha)$, $h^{(j)}(\alpha) = r_{j,\alpha}$.

Isomorphisms among the $R_{L,\gamma,\delta}$ are easily determined. For real $c,d$ with $c \neq 0$, let $\psi_{c,d}$ be the ring automorphism determined by $x \mapsto cx+d$. The ring automorphisms of  $L[x]$ are exactly those of the form $\psi_{c,d}$ (restricted to $L[x]$) where $c$ and $d$ belong to $L$. These induce ordered ring isomorphisms $R_{L,\gamma,\delta} \to R_{L,\gamma',\delta'}$ where $\psi_{c,d}$ sends the two-element set $\brcs{\gamma, \delta}$ to $\brcs{\gamma', \delta'}$. Since the pure traces are determined, even topologically, by the ordering, any ordered ring isomorphism will have to have this property, so that the order-isomorphism classes (for $L$ fixed) are precisely the orbits of $\brcs{\gamma,\delta}$ under $\psi_{c,d}$ where $c,d \in L$ (and $c\neq 0$). 

More interesting is what happens when we change the base field. Let $L \subset K \subseteq \R$ be a proper field extension of $L$, inside $K$. We form the tensor product, $R_{L,\gamma,\delta} \otimes_L K$ as $L$-vector spaces; of course, as an $L$-algebra, this is just $K[x]$. We can impose the tensor product ordering, by considering the cone generated by the pure tensors of the form $f \otimes k$ where $f \in R_{L,\gamma,\delta} $ and $k \in K^+$ (ordering inherited from $\R$) [we would have to check that this is a proper cone, as in [GH], but since $\R$ is an injective $K$-module---as $K$ is a field---the same argument works.] Although $R_{L,\gamma,\delta} $ is by definition archimedean (for ordered abelian groups with order unit, this is equivalent to $\tau(g) \geq 0$ for all pure traces $\tau$ implies $g \geq 0$), the tensor product is not. 

We simply note that the zero sets of positive elements of the tensor product are finite subsets of $[\gamma,\delta]$ that consist of $L$-algebraic numbers and which are relatively closed under conjugacy (that is, if $\alpha$ and $\beta$ are algebraic conjugates and both are in the interval, then $\alpha$ is in a zero set of a positive element if and only if $\beta$ is in the same one). In particular, if a singleton $\brcs{\alpha}$ arises as a zero set of a positive element, then $\alpha \in L$. So select $\eta$ in $K \setminus L$; since $L\eta \setminus \brcs{0}$ is dense in the reals, we can find $\alpha = a\eta$ in $[\gamma,\delta]$ for some nonzero $a$ in $L$. The element 
$(x- \alpha)^2$ (strictly speaking, $x^2 \otimes 1 - 2 ax \otimes \eta + a^2 \otimes \eta^2$) is in the tensor product. It is nonnegative as a function on the interval, hence on the pure traces; but it cannot be in the positive cone of the tensor product as its zero set is the singleton $\brcs{\alpha}$, which is not closed under conjugation \wrt $L$.

What can be shown is that if $L_i \subset K_i \subseteq \R$ are two proper field extensions, then the ordered tensor products, $R_{L_i,\gamma_i,\delta_i} \otimes_{L_i} K_i $ are isomorphic (as ordered rings) if and only if there is a field isomorphism $\Arrow \sigma; K_1 . K_2$ whose restriction to $L_1$ yields an isomorphism to $L_2$,  for which there is a $\phi_{c,d}$ with $c,d \in L_2$   compatible with the maps, mapping the endpoints to endpoints. 

Similar results apply to the trigonometric polynomial ring (another of Fuchs' examples of ordered rings satisfying interpolation)---the criterion in the technical lemma applies even more generally to real analytic functions.

\comment

\vskip 4pt \noindent {\bf Appendix to the Appendix} Tensor products of modules over partially ordered rings

\vskip 2pt \noindent Let $A$ be a partially ordered commutative ring for which $1$ is an order unit. Let $M$ and $N$ be ordered $A$-modules, that is, they are partially ordered abelian groups that are also (unital) $A$-modules, and satisfy $M^+ A^+ \subseteq M^+$. We may thus form the $A$-module tensor product, $M \otimes_A N$ which will be abbreviated $M \otimes N$. This is an $A$-module, and we may put a candidate positive cone on it, the {\it tensor post-ordering,} the set of elements of the form
$\sum a_i m_i \otimes n_i$ where $a_i \in A^+$, $m_i \in M^+$, and $n_i \in N^+$. All the properties of being an ordered $A$-module hold, except possibly that the cone is proper.

When $A = \Z$ (considering ordered abelian groups), it was shown in [GH] that the cone is proper. The tools used were the existence of nontrivial traces, and the fact that $\R$ is injective as an abelian group. We can show that even in the general case, the analogous results can be obtained (with more effort), so as to adapt the proof in [GH]. 

\Lem Proposition. If $A$ is a commutative partially ordered ring with $1$ as an order unit, and $M$ and $N$ are partially ordered $A$-modules, then the $A$-module tensor product $M \otimes_A N$, equipped with the tensor product post-ordering, is a partially ordered $A$-module. 

The following is an easy adaptation (using a standard technique) of a result about pure traces on order ideals.

\Lem Lemma. Let $A$ be a partially ordered ring with $1$ as order unit, and let $(M,u)$ be an ordered $A$-module with order unit. If $\tau$ is a pure trace on $(M,u)$, then there exists a pure trace, $L$, of $A$ \st for all $m$ in $M$ and $a$ in $A$, we have $\tau (ma) = \tau(m) L(a)$.

\Pf Let For each $b$ in $A^+$, let $\Arrow \phi_b; M.M$ be right multiplication by $b$; it is clearly an order-preserving $A$-module endomorphism. Form the new trace (which could be zero), $\tau \circ \phi_b$. From $1$ being an order unit for $A$, that there exists positive integer $K$ \st $b \leq K1$, and then it follows that as positive additive functions on $M$, the $\tau \circ \phi_b \leq K \tau$. From $\tau$ being pure, there exists a unique nonnegative real number $\gamma(b)$ (which could be zero) \st $\tau \circ \phi_b = \gamma (b) \tau$. It is easy to check that on $A^+$, $\gamma$ is additive, multiplicative, and sends $1$ to the real number $1$. Moreover, if $b_1 - b_2 = b_3 - b_4$ (with $b_i$
 in $A^+$), then $b_1 + b_4 = b_2 + b_3$, so $\gamma(b_1) + \gamma(b_4) = \gamma(b_2) + \gamma(b_3)$, whence $\gamma(b_1) - \gamma(b_2) = \gamma(b_3) - \gamma(b_4)$, and thus $\gamma$ extends uniquely to a  trace $L$, on $A$, and it is clearly multiplicative, hence pure;  moreover, $\tau (ma) = \tau \circ \phi_a (m) = \tau (m) L(a)$ if $a$ is in $A^+$, and it is trivial that this property holds for all $a$ in $A$. \qed

The pure traces on $(A,1)$ are exactly the multiplicative ones. We can view $\R$ as an $A$-module (structure induced by $L$), by setting $a\cdot r$ to be $L(a)r$. Let $S\subseteq \R$ be the image of $A$ under $L$. Now we forget about the ordering.

\Lem Lemma. Let $S$ be a unital subring of $\R$, and endow the latter with the obvious $S$-module structure, to create the module $\R_S$. Then $\R_S$ is injective (as an $S$-module). 

\Pf Let $P \subset Q$ be $S$-modules, and suppose there is an $A$-module homomorphism $\Arrow \phi; P.\R$. We wish to extend $\phi$ to $Q$. Let $T(B)$ be the $S$-torsion  submodule of the module $B$, that is, $\Set{b\in B}{\exists s \in S\setminus \brcs{0}, bs = 0}$, and let $\overline B$ denote the torsion-free module $B/T(M)$. It is obvious (since $S$ is a domain) that $T(P) \subseteq \ker \phi$, and moreover, $T(Q) \cap P = T(P)$. Hence the induced map $\overline P \to \overline Q$ is an embedding of $S$-module, and of course, $\phi$ induces a now one to one $S$-module map $ \overline P \to \overline Q$. Let $K$ be the field of fractions of $S$ inside $\R$. We can enlarge $\overline P$ into a $K$-vector space (in such a way that the map to the reals extends), by noticing that every nonzero $s$ in $S$ acts in a now one to one fashion on $\overline P$ (and also on $\overline Q$). 

We thus take the direct limit over all maps $\times s$ with each $s$ in $S$ repeated infinitely many times, doing this both to $\overline P$ and $\overline Q$ and the map between them. The result (which can be identified with $\overline P \otimes_S K \to \overline Q \otimes_S K$) is an inclusion of $K$-vector spaces, $P_K \to Q_K$. Moreover, $\overline P $, $\overline Q $ sit inside their counterparts in such a way that the restriction of the map to $\overline P$ yields the map $ \overline P \to \overline Q$. 

For every $s$ in $S$, the map $\overline \phi$ extends by simply defining $ps^{-k} \mapsto \overline\phi(p)/s^k$ (on the right side, viewing $1/s^k$ as a real number). Thus $\overline \phi$ extends to a map $P_K \to \R$. This is now a map as $K$-vector spaces, and since everything is injective over a field, we can extend this map to $Q_K $, and thus by restriction to $\overline Q$, and following the arrows yields the desired extension. 
\qed

final comment
\endcomment

\comment
Suppose $f_i$ and $g_j$ ($i,j = 1,2$) are four real polynomials \st for all $i,j$, $f_i \leq g_j$ as functions on the unit interval. Let $F = f_1 \vee f_2$ and $G = g_1 \wedge g_2$, so as continuous functions on the interval, $F \leq G$. We may assume that neither $f_i$ is equal to either of the $g_j$. Then the zero set of function $G-F$ must consist of a finite set of points $Z = \brcs{x_1, \dots, x_k}$ (in other words, $Z$ consists of the points in the unit interval where at least three of the polynomials $\brcs{f_1, f_2, g_1, g_2}$ agree. 

Suppose that $f_2(x_l) = f_1 (x_l) = g_1(x_l) \leq  g_2 (x_l)$ for some $l$ (it may also happen that all four functions agree at $x_l$). If we consider the functions $g_1 - f_1$ and $g_1-f_2$ (and, if necessary, the other two differences), we note that these are nonnegative polynomials vanishing at $x_l$. In particular, if we control enough derivatives of a polynomial $h$ at $x_l$,  with $h (x_l) = f_1(x)$, then automatically it will follow that there will be a neighbourhood of $x_l$ (in $[0,1]$) \st  $h$ is squeezed between $f_1$, $f_2$ and $g_1$ (on the neighbourhood); it could be a one-sided neighbourhood if $x_l $ is $0$ or $1$. By Hermite interpolation, there exists a polynomial with specified values and values of derivatives at each of the $x_l$. 
Hence there exists a polynomial $h$, together with an open set $U$ containing $Z$, such that for all $u$ in $U$, $F(u) \leq h(u) \leq G(u)$. [Obviously the same process works if instead $g_2(x_l) = g_1 (x_l) = f_1(x_l) \geq f_2(x_l)$ for some $l$.]

For each $x_l$, we have that the four differences $h- f_i$ and $g_j - h$ are nonnegative on a neighbourhood of $x_l$ and some of them vanish; let $m_l $ be the maximum order of the zero at $x_l$ of the four functions; since they are nonnegative, $m_l$ is even, except possibly if $x_l$ is zero or one. Let $M $ be whichever of  $\max\brcs{m_l}_{l=1}^k + 2$ or $\max\brcs{m_l}_{l=1}^k+3$ is even.

Set $r = \prod_{l=1}^k (x-x_l)^M$, so that $r$   is nonnegative on the unit interval, and vanishes only on $Z$. Form the rational functions $(g_j -h)/r$ and $(f_i - h)/r$. Each has pole of order at least $2$ at each of the $x_l$; we notice that the $(g_j -h)/r$ are nonnegative on $U \setminus Z$, so approach infinity as points approach an $x_l$; on the other hand, $(f_i - h)/r$ are nonpositive on $U \setminus Z$, so approach $-\infty$ at points approaching the $x_l$. Thus, for each $l$, there exists a punctured neighbourhood of $x_l$, $W_l$, and a positive integer $N_l$   \st $(g_j -h)/r - (f_i - h)/r > N_l$ on $W_l$; set $N = \inf \brcs{N_l}$ and $W = \cup W_l$. 

Each of $(f_j - h)/r$ is bounded above, say by $A > 0$, i.e., $A > (f_j - h)/r$ on $[0,1]\setminus Z$, and similarly, $(g_i -h)/r$ is bounded below by $-B$ (where $B$ is a positive real number). We may assume $A = B$. Form the functions
$F_i = -A\vee (f_i - h)/r$ and $G_j = A\wedge (g_j - h)/r$. Then $F_i \leq G_j$, each of the four is continuous, and moreover, the inequalities are pointwise strict (that is, there exists $\delta > 0$ \st for all $t$ in $[0,1]$, $G_j (t) - F_i (t) \geq \delta$). By Weierstrass' approximation theorem, we can uniformly approximate the continuous function $(G_1 \wedge G_2)/2 + (F_1 \vee F_2)/2$ by a polynomial, $s$, to within $\delta/2$. It follows immediately that $F_i \leq s \leq G_j$, and all the more so, the pointwise inequalities (on $[0,1] \setminus Z$)
$$
\frac{f_i - h}r \leq s \leq \frac{g_j - h}r.
$$
Thus $f_i \leq h+rs \leq g_j$ pointwise on $[0,1] \setminus Z$; by continuity, these inequalities are valid on all of $[0,1]$. \qed

The use of Hermite interpolation in Fuchs' outline proof is slight overkill (while an easy theorem to prove, it usually comes with additional restrictions, such as on degree; this is unnecessary here, and in fact the Chinese remainder theorem is all we really need). 

Let $f_1,f_2,g_1,g_2 \in \R[x]$ \st $g_i \leq f_j$ on the unit interval. Suppose that no $g_i$ equals an $f_j$; hence the set of points, $Z$, consisting of $t$ in the interval where $\sup \brcs{f_1, f_2} (t) =\inf \brcs{g_1, g_2} (t)$ is  finite. Suppose $\alpha$ is in $Z$ and is an interior point of the unit interval.

We will show there exists an integer $m \cong m(\alpha)$ to together with a nonempty  open interval $(a(\alpha), b(\alpha))$ and a bunch of real numbers, $r_j$ \st if $h$ is any real polynomial satisfying
 
\item{(i)} $h(\alpha) = \inf \brcs{g_1, g_2} (\alpha)$
\item{(ii)} $h^{(j)}(\alpha) = r_j$ for $1 \leq j \leq m$
\item{(ii)} $h^{(m+1)}(\alpha)  \in (a(\alpha), b(\alpha))$,

\noindent then there exists a neighbourhood $U \equiv U(\alpha,h)$ \st restricted to $U$, we have $f_i \leq h \leq g_j$.

Subtracting $f_2$ from each of the terms, we may reduce to the case that $f_2$ is identically zero.

Select $\alpha$ in $Z \cap (0,1)$. We may relabel the $g$s so that $g_2(\alpha) \geq g_1(\alpha) = f_1(\alpha) $.

\noindent Case 1. {\it Suppose   that $\inf  \brcs{g_1, g_2} (\alpha)\neq 0$, hence is greater than zero.} Since $g_1 - f_1 \geq 0$ and vanishes at $\alpha$, we may factor it as $(x-\alpha)^{2m_1}L_1$ where $L_1 (\alpha) > 0$, for some positive integer $m_1$ (the zero must be of even order, since the function is nonnegative in a neighbourhood of $\alpha$); if in addition, $g_2(\alpha) = g_1(\alpha)$, we have a corresponding factorization
$g_2 - f_1 = (x-\alpha)^{2m_2}L_2$. If $m_1 > m_2$, or if $m_2$ does not exist (that is, is zero), set $m = 2m_1-1$, and set $r_j = f_1^{(j)}(\alpha)$. Then $h - f_1 = 0 + c(x-a)^{2m_1} + \Oh{(x-a)^{2m_1 +1}}$ on a neighbourhood of $\alpha$, so that if $c> 0$, then $h  \geq f_1,0  $ on a neighbourhood of $\alpha$; and moreover $g_i -h  = g_i - f_1  - h + f_1 = (x-\alpha)^{2m_i}L_i - c(x-a)^{2m_1} $, and this will be nonnegative on a neighbourhood of $\alpha$ if $c < L_1 (\alpha)$.

If $m_1 = m_2$, define the $r_j$ as before; then the condition is simply that $c < \min\brcs{L_1 (\alpha), L_2 (\alpha)}$.
And if $m_1 < m_2$, set $m = 2m_2 -1$ and $r_j = f_1^{(j)} (\alpha)$, as in the case $m_2 < m_1$.

\noindent   Case 2. {\it Suppose that $g_1 (\alpha) = 0$}. Since $g_1 \geq 0$, we may also write $g_1 = (x-\alpha)^{2n_1}M_1$ where $M_1 (\alpha) > 0$, and similarly, if $g_2 (\alpha)$ is also zero, then $g_2 =  (x-\alpha)^{2n_2}M_2$. Also, $f_1 (\alpha) \leq 0$, so either it is strictly negative, or it is also zero there. In the latter case, we have $f_1 = (x-\alpha)^t P$ where $P(\alpha) \neq 0$.

{\parindent = 5em
\item{Case 2a} {\it Suppose that $f_1(\alpha) < 0$}. Set $m = 2\max{n_1,n_2}-1$ and $r_j = 0$. If $m_i > m_j$, $0 < c < L_i(\alpha_i)$ will do, whereas if $n_1 = n_2$, then $0 < c < \min{L_1(\alpha),L_2(\alpha)}$ is sufficient.

\item{Case 2b} {\it Suppose that $f_1 (\alpha) = 0$.} We may write $f_1 = (x-\alpha)^t P$ where $P(\alpha) \neq 0$. We note that if $n_1 \neq n_2$, then there is a neighbourhood of $\alpha$ on which $g_1$ is comparable to $g_2$, i.e., $g_1 \leq g_2$ (after a possible relabelling) on this neighbourhood, so we can disregard $g_2$.
\itemitem{Case 2bi} {\it Suppose that $t$ is even.} If  $P(\alpha) < 0$, proceed as in Case 2a. Otherwise, we must have $2n_i \leq t$.  Divide everything in sight by $(x-a)^{2n_1}$, arriving at a situation where the value at $\alpha$ is not zero---this is case 1, and we simply modify the $m$ by adding $t/2$ to it, and the initial segment of the $r_j$ are all set to zero. If $n_1 = n_2$, we can similarly divide by $(x-a)^{2n_1}$, and again reduce to case\, 1.

 \itemitem{Case 2bi} {\it Suppose that $t$ is odd.} Then $f_1$ changes sign at $\alpha$, but we can still work on  the side on which $f_1$ is positive. It follows from $g_1 \geq f_1$ that $2n_1 \leq t$, hence $2n_1 < t$. Set $m = 2n_1 -1$ and $r_j =0$; the only constraint is that $0 < c < M_1{\alpha}$ (the $n_2$ portion can be ignored---if $n_1 = n_2$, we may have assumed that $M_1(\alpha) \leq M_2(\alpha)  $, and if $n_2 < n_1$, $g_2$ is locally greater than $g_1$ anyway.
\qed
}

Now it is obvious that the map $R = \R [x] \to \R^{m+2}$ given by $f \mapsto (f^{(j)}(\alpha))_{j = 0}^{m+1}$ is a vector space surjection (the kernel consists of all polynomials whose derivatives vanish at $\alpha$, that is, the ideal $I= (x-a)^{m+2}$, and $R/I$ is equidimensional with the right side). Slightly more interesting is what happens over the rational numbers. 

Set $R = \Q[x]$. Suppose that $\alpha$ is a algebraic  (satisfies a nonconstant polynomial equation with rational coefficients), and let $p$ be the monic irreducible polynomial it satisfies, and set $K = \Q[\alpha] \iso R/pR$. The map 
$R = \Q[x] \to K^{m+2}$ defined in the same way ($f \mapsto (f^{(j)}(\alpha))_{j = 0}^{m+1}$) is a surjection of $\Q$ vector spaces (as before, the kernel consists of the ideal $I = p^{m+1}R$, and $R/I$ has dimension $(m+2)\dim_{\Q} K$, which is equidimensional with the right side. This simply says that we given a bunch of elements of $K$, $\brcs{r_0, r_1, \dots r_{m+1}}$, there exists a rational polynomial $q$ \st $q^{(j)} (\alpha) = r_j$ for $ j=0,1,\dots,m+1$. 

Now we use the Chinese remainder theorem: if $I_{s}$ is a finite family of ideals of a ring $R$ \st that for all $s \neq t$, we have $I_s + I_t = R$, then  the natural map $R \to \prod_s R/I_s$ is onto (i.e., the the quotient $R/I$ is naturally isomorphic to the product of the quotients).

If $R = \R [x]$, let $\alpha_s$ be a finite set of distinct real numbers, and for each $s$, let $r_{js}$ be a finite initial segment of real numbers ($0 \leq j \leq N(s)$). Then there exists a polynomial $f$ \st $f^{(j)} (\alpha_s) = r_{sj}$.

Let $p_s = (x-\alpha_s)^{N(s}+1}$, and set $I_s$ to be the principal ideal generated by $p_s$. Then the ideals are comaximal, hence there given the elements $r_s$, there exists $f$ in $\R[x]$ \st $f - r_s \in I_s$, i.e., $f - r_s = p_s h $ for some polynomial $h$. Since $D^i p_s (\alpha_s) = 0$ for $0 \leq i \leq N(s)$, the derivatives of $f$ match those of $r_s$ at $\alpha_s$. We can select $r_s$ to be a polynomial whose derivatives evaluated at $\alpha_s$ are arbitrary, as we noticed above.

In the case of the rational polynomial ring, that is, $R = \Q[x]$, we have to be more careful. First, the points $\alpha_s$ must be algebraic---evaluating a polynomial at a transcendental determines the polynomial completely (assuming the transcendental is known), and moreover the specified values must lie in $K_{\alpha} = \Q{\alpha_s}$. Second, set  $K_s = \Q[\alpha_s]$, and let $\sigma$ be an element of the Galois group of the algebraic closure that maps $K_s$ to the reals, say $\sigma (\alpha_s) = \beta \in \R$. If $\beta = \alpha_t$ (for some  $t \neq s$), then $f^{(j)}(\alpha_t) = \sigma (f^{(j)}(\alpha_s))$---in other words, the values of the derivatives at $\alpha_t$ are determined by those at $\alpha_s$. Say $\beta$ is an (algebraic) conjugate to $\alpha$ if such a $\sigma$ exists, i.e., their irreducible monic polynomials are the same. 
Then we have the following; necessity of the conditions is obvious (except in the trivial case when $|S| =1$ and $N(s) =0$). 

 Suppose that $\brcs{\alpha_s}_{s\in S}$ is a finite set of real numbers, together with data $\brcs{r_{sj}}_{j=0}^{N(s)}$ \st the following  conditions hold:
\item{(a)} each $\alpha_s$ is algebraic;
\item{(b)} for all $j$ and $s$, $\r_{sj} \in K_s = \Q[\alpha_s]$;
\item{(c)} if $\alpha_s$ is an algebraic conjugate of $\alpha_s$, say with $\sigma (\alpha_t) = \alpha_s$, then $N(s) = N(t)$ and $r_{js} = \sigma({r_{jt}})$ for all relevant $j$.

\noindent Then there exists $f$ in $\Q[x]$ \st $D^j f(\alpha_s) = r_{js}$ for all relevant $j$ and $s$.

\Pf For any specific $t$, there exists $f_t$ in $\Q[x]$ whose derivatives match $r_{jt}$ by (b) and the ontoness of the map $R:= \Q[x] \to K_t^{N(t)+1}$. If $\sigma(\alpha_t) = \alpha_s$, applying $\sigma$ to each rational polynomial $f_t^{(j)}$ evaluated at $\alpha_s$, we have $f_t^{(j)}(\alpha_s) = f_t^{(j)}(\sigma(\alpha_t)) = \sigma(f_t^{(j)}(\alpha_t))$, and the last equals $r_{js}$ by (c). Hence $f_s$ automatically serves as an interpolant at $\alpha_t$. 

Let $T$ be a subset of $S$ \st that $\brcs{\alpha_s}_{s\in T}$ contains no pairs of conjugate numbers, but every element of $S$ is conjugate to an element of $T$. For each $t$ in $T$, let $p_t$ be the monic  polynomial of $\alpha_t$ irreducible  (over the rationals), and let $I_t = p_t^{N(t) +1}R $; then the Chinese remainder theorem applies, and thus there exists a polynomial $f$ in $R$ \st $f - f_t \in I_t$ for each $t$ in $T$. Now the values of the derivatives match, since $D^{N(t)}$ kills $I_t$. Since every $\alpha_s$ is conjugate to  $\alpha_t$ for unique $t$ in $T$, the previous paragraph yields that $f - f_s$ vanishes at all the relevant derivatives evaluated at $\alpha_s$
\qed
 
We have one more (wellknown) step to deal with the rational polynomial case. Let $K$ be a formally real finite dimensional extension field of $\Q$, and let $\brcs{\sigma}_{\Sigma}$ be a  family of homomorphisms $K \to \R$ (at least one exists by formal reality). Then the image of $K$ under the map $K \to \R^{\Sigma}$ is dense. We note that $\Sigma$ is automatically linearly independent, so both $|\Sigma| <\infty$ and density follow immediately. 

Suppose that $f_1, g_1, g_2$ are in $\Q[x]$ and $0, f_1 \leq g_1, g_2$ on an interval of the form $[a,b]$ (unlike the situation over the reals, it actually makes a difference which endpoints are chosen). Let $S$ be the set of interior points at which $g_1 \wedge g_2$ agrees with $ f_1 \vee 0 $. Let $T$ be a subset containing exactly one of the algebraic conjugacy classes of the $\alpha_s$, and partition $S $ into subsets consisting of single conjugacy classes, i.e., write $S = \cup_{t \in T} S_t$ where $ s,s' \in S_t$ if and only if $\alpha_s$ and $\alpha_{s'}$ are conjugate. Fix $t$, and for each $s \in S_t$, we may find an element $u_s$ in $K_s$ \st under the corresponding real embedding, $a(\alpha_s) < u < b(\alpha_s)$. Set $K_t = \Q[\alpha_t]$ ($\brcs{t} = T \cap S_t$) and let $\Sigma$ be the set of homomorphisms corresponding the conjugates, $\alpha_s$ ($s\in S_t$). By density, we may find $k_t \in K_t$ \st the image of $k$ in $K_s \subseteq \R$ lies in the interval $(a(\alpha_s), b(\alpha_s))$. Now define for each $t$, $r_{jt}$ from the conditions determined by the polynomials, and a final one (for the last derivative), $r_{N(t)+1,t} = k_t$. We still have to show that $\sigma(r_j) $ is what it is supposed to be. ...
\endcomment
\comment
If $2n_i = t$, we obtain either $M_i (\alpha) = P(\alpha)$ and $m_i > n_i$ or $M_i (\alpha) > P(\alpha)$ and $m_i = n_i$. In the 

If $$

If $t$ is even and $P(\alpha) > 0$, set $m = 2\max\brcs{m_1,m_2} -1$ and $r_j =0$. If now $m+1 = t$, the constraint is simply $P(\alpha)$

}

 If $t$ is even and $P(\alpha) < 0$, then $f_1$ is nonpositive on a neighbourhood of $\alpha$, and we can throw it away, if we ensure that $h$ is nonnegative on a neighbourhood. To this end, $m = 2n_1 - 1$ if $n_2$ does not exist, and $2\max{2n_1,n_2} -1 $ if it does, and set $r_j = 0$. Then $h = c(x-a)^{\max{2n_1,2n_2} + \Oh{(x-a)^{\max{2n_1,2n_2 +1}}$, so if $c > 0$, $h \geq 0$ on a neighbourhood, and if $c < M_i (\alpha)$, then $h \leq g_i$ on a neighbourhood.

If $t$ is even and $P(\alpha) > 0$, then $f_1$ is nonnegative on a neighbourhood, and $g_i \geq f_1$ entails $n_i \leq t/2$; select $m = t-1$, and verification is simple.

Finally, if $t$ is odd, then $f_1$ changes sign at  $\alpha$; in particular, on one side of $\alpha$, it behaves like $|(x-a)|^t |P(\alpha)| + \Oh{(x-a)^{t+1}}$; since $g_i \geq f_1$ everywhere, it follows that $t \geq 2n_1, 2n_2$.

 of the four functions must agree at $\alpha$; for convenience, assume $g_1,g_2$, and $f_1$ do (otherwise, we can just multiply by $-1$ and reverse the implications), and then $g_2(\alpha) \geq g_1(\alpha) $ (and it is possible that all four agree. Consider $g_i - f_1$; as it vanishes at $\alpha$ and is nonnegative (everywhere in the interval), in particular on a neighbourhood of $\alpha$, we may factor it as $(x-\alpha)^{2m_i}L_i$ where $L_i$ is a polynomial \st $L_i (\alpha) > 0$. If $g_2 (\alpha) = g_1(\alpha)$, we also have factorizations, $g_2 - f_i = (x-\alpha)^{2n_i}M_i$. Set $m(1) $ to be the maximum   of $m_1, m_2$ together with $n_i$ (if the latter exist).  

If $m(1)$ equals at most one of the $m_i$ (and the $n_i$, if they exist), then $x-\alpha)^{2m(1)} $ 
\endcomment

\subtitle References

\long\def\Reff[#1] #2, #3, #4\par{\vskip 1pt \item{[#1]} #2, {\it #3,}
#4\par} {\parindent = 2.5em %\everypar={\Reff } 

\Reff [EHS] Edward G Effros{, David E Handelman,} and Chao-Liang Shen,
Dimension Groups and Their Affine Representations, 
American Journal of Mathematics 102 (1980) 385--407.

\Reff [F1] L Fuchs, Riesz groups, Ann Scuola Norm Sup Pisa (3) 19  (1965) 1--34.

\Reff [F2] L Fuchs, Riesz rings, Mathematische Annalen 166 (1966)  24--33.

\Reff [F3] L Fuchs, Riesz vector spaces and Riesz algebras,  Queen's papers in pure and applies mathematics, No\,1, Queen's University, Kingston  (1966).

\Reff [G] K Goodearl, Partially ordered abelian groups with interpolation, Mathematical series and monographs \#20, American Mathematical Society (1986).

\Reff [GH] K Goodearl and D
Handelman, Tensor products of dimension groups and K${}_0$ of unit regular rings, Canad J Math 38 (1986) 633--658.

}  
 
\vskip 6pt \noindent Mathematics Department, University of Ottawa, Ottawa
ON K1N 6N5, Canada; dehsg\@uottawa.ca

\end